\documentclass[a4paper,12pt,reqno]{article}

\usepackage{t1enc,amssymb,amsbsy,amsthm,amsmath,amsfonts,graphics}
\usepackage[normalem,ULforem]{ulem}
\usepackage{color}
\newtheorem{lemma}{Lemma}
\newtheorem{theorem}{Theorem}
\newtheorem{corollary}{Corollary}

\usepackage{amsfonts,amsmath,amssymb,bbm}
\usepackage{setspace} 
\def\proof{\noindent{\bf Proof:}\hskip10pt}        
\def\QED{\hfill $\Box$}

\font\tenmath=msbm10 scaled 1200
\font\sevenmath=msbm7 scaled 1200
\font\Fivemath=msbm5 scaled 1200
\newfam\mathfam \textfont\mathfam=\tenmath
\scriptfont\mathfam=\sevenmath \scriptscriptfont\mathfam=\Fivemath

\def\R{\mathbb{R}}
\def \1{1 \mkern -6mu 1} 

\def\E{\mathbb{E}}
\def\P{\mathbb{P}}

\def\R{\mathbb{R}}

\newcommand{\dis}{\displaystyle}

\newcommand{\tg}{t\geq0}

\newcommand{\lf}{\left(}
\newcommand{\ri}{\right)}

\newcommand{\lff}{\left[}
\newcommand{\rii}{\right]}

\newcommand{\RR}{\mathbb{R}}

\def\R{\mathbb{R}}
\def \1{1 \mkern -6mu 1} 

\def\E{\mathbb{E}}
\def\P{\mathbb{P}}

\def\R{\mathbb{R}}

\newcommand{\law}{\overset{\textrm{(law)}}{=}}

\newcommand{\deff}{\overset{\textrm{def}}{=}}

\begin{document}

\title{\textbf{Some two-dimensional extensions of Bougerol's identity in law for
 the exponential functional of linear Brownian motion}}

\author{J. Bertoin\thanks{Institut f\"ur Mathematik, 
Universit\"at Z\"urich, 
Winterthurerstrasse 190, 
CH-8057 Z\"urich, Switzerland. \eject  E-mail:
{\tt jean.bertoin@math.uzh.ch}} \and 
D. Dufresne\thanks{Centre of Actuarial Studies, 
University of Melbourne, 3010 Victoria,
Australia. \eject E-mail:{\tt dufresne@unimelb.edu.au}} \and 
M. Yor\thanks{Laboratoire de Probabilit\'es et Mod\`eles Al\'eatoires, 
UPMC, 4 place Jussieu, F-75252 Paris cedex 05, France. E-mail:
{\tt deaproba@proba.jussieu.fr}} \thanks{Institut Universitaire de France}}

\maketitle

\begin{abstract}  We present a two-dimensional extension of an identity in distribution due to Bougerol \cite{Bou} that involves the exponential functional of a linear Brownian motion. Even though this identity does not extend at the level of processes, we point at further striking relations in this direction.

\end{abstract}

{\bf Key words:} Brownian motion, exponential functional, Bougerol's identity, local time, Bessel processes. 

\section{Introduction}\label{intro}

\textbf{(1.1)} To a linear Brownian motion $(B_s,\,s\ge0)$ starting from 0, we associate the exponential functional
\begin{equation*}
    A_t=\int_0^tds\exp(2B_s)\,,\qquad t\ge0.
\end{equation*}
The distribution of $A_t$ is made accessible thanks to Bougerol's identity in law
\begin{equation}\label{1}
    \textrm{for fixed}\quad t\;,\quad\sinh(B_t)\law\beta(A_t),
\end{equation}
where $(\beta(u),\,u\ge0)$ denotes a Brownian motion which is independent of $(B_s\,,s\ge0)$, hence of $A_t$.
Assuming \eqref{1},  elementary computations yield the following
characterization of the law of $A_t$:
\begin{equation}\label{2}
    \E\lff\frac{1}{\sqrt{A_t}}\exp\lf-\frac{x^2}{2A_t}\ri\rii=\frac{a'(x)}{\sqrt{t}}\exp\lf-\frac{a^2(x)}{2t}\ri\,,\qquad x\in\RR\,,
\end{equation}
where
\begin{equation*}
  a(x)=\arg \sinh (x)\equiv\log(x+\sqrt{1+x^2}) \quad \hbox{and} \quad  a'(x)=\frac{1}{\sqrt{1+x^2}}\,.
\end{equation*}

For further reference, we note some simple, but useful, consequences of (\ref{2}), {\it i.e.}:
\begin{equation}\label{3}
 \E\left({1\over\sqrt{A_t}}\right) = {1\over\sqrt{t}}\,,
 \end{equation}
and, differentiating both sides with respect to $t$:
\begin{equation}\label{4}
    \E\left(\frac{\exp(B_t)}{A_t^{3/2}}\right)\overset{(\ast)}{=}\E\left(\frac{\exp(2B_t)}{A_t^{3/2}}\right)={1\over t^{3/2}}\,,
\end{equation}
where $(\ast)$ is obtained by time reversal of $(B_s,s\leq t)$
from time $t$.

\noindent\textbf{(1.2)}  It took some time, despite the original proof in \cite{Bou}, to understand simply and deeply why (1) holds. In \cite{ADY}, one finds the following arguments, among which the (essential) time reversal one:
\begin{align*}
\text{for fixed }t\,,\  \beta(A_t)\ \ \text{is distributed as:}\ \ \int_0^t\exp(B_s)\,d\beta(s)\,,
\end{align*}
which, by time reversal (at time $t$) is also distributed as
\begin{align}\label{timereversal}
\exp(B_t)\int_0^t\exp(-B_s)\,d\beta(s)\,.
\end{align}
Now, it is easily shown, using It\^o's formula, that  the {\em process} in \eqref{timereversal} is distributed as the {\em process} $(\sinh(B_t)\,,\, t\ge0)$.

\noindent{\bf(1.3)} In the present paper, we obtain an extension of \eqref{1}, by considering the two-dimensional vector $(\sinh(B_t),\sinh(L_t))$, where $(L_t\,,\,t\ge0)$ denotes the local time at 0 of $B$. Our main result is:

\begin{theorem}\label{theo1}
For fixed $t$, the three following  two-dimensional random variables
are identically distributed:
\begin{align}\label{ids}
  (\sinh(B_t)\,,\,\sinh(L_t))\ \law \ (\beta(A_t)\,,\,\exp(-B_t)\lambda(A_t))\ \law\ (\exp(-B_t)\beta(A_t)\,,\,\lambda(A_t))\,,
\end{align}
where $(\beta(u),u\geq0)$ is a one-dimensional Brownian motion, with local time
at 0, $(\lambda(u),u\geq0)$, and $\beta$ is independent from $B$.
\end{theorem}
 It may be interesting to observe right now that Tanaka's formula shows that the local time at level $0$ and time $t$ of the process $(\sinh(B_s), s\geq 0)$ is simply $L_t$, whereas that of the process $(\exp(-B_s)\beta(A_s), s\geq 0)$ can be expressed as 
$\int_0^t\exp(-B_s) d\lambda(A_s)$. Hence we have also the identity in distribution between two-dimensional processes
\begin{equation}\label{eqproc}
\left(\sinh(B_t), L_t
\right)_{t\geq 0} \law \left(\exp(-B_t)\beta(A_t), \int_0^t\exp(-B_s) d\lambda(A_s)\right)_{t\geq 0}\,.
\end{equation}
We stress that  \eqref{ids} cannot be extended to the level of processes; see the forthcoming Section 2.2. Hence the two identities  in distribution \eqref{ids} and \eqref{eqproc} differ profoundly.

Theorem \ref{theo1} is proved in Section 3.
In Section 2, we discuss a number of consequences and equivalent statements to that of Theorem 1. For instance, the well-known equivalence in law, due to Paul L\'evy, between the processes $((\overline B_t-B_t,\overline B_t)\,,t\ge0)$ and $((|B_t|\,,L_t)\,,t\ge0)$ allows to present a version involving the supremum $\overline B_t=\sup_{s\le t}B_s$ instead of the local time version of Theorem \ref{theo1}.

Apart from this, Section 2 consists in the statements and discussions of  four other theorems. Roughly speaking, these theorems were motivated by our desire to understand whether in \eqref{ids}, the two extreme identities hold for processes. This question has now been solved in the negative (see \cite{BY}), but nonetheless there are some rather remarkable identities between jump intensity measures, which are described in Theorems 2-5, and which made us believe for some time in a $2$-dimensional process identity extending \eqref{ids}. We let the reader discover the precise statements of these theorems in Section 2;
their proofs are found in Section 4.

\noindent{\bf (1.4)} Before we get into the precise arguments of the proofs of our various theorems, we should like to present an overall appreciation of the present work, by making a parallel with the way our understanding of Bougerol's identity \eqref{1} has improved: in \cite{Yb}, a Mellin transform proof was given, based on the identity in law \eqref{eq12} below. Later, in \cite{ADY},  a time-reversal argument and stochastic calculus proof of \eqref{1} were found. We estimate that, at the moment, our understanding of Theorem \ref{theo1} lies at the level of \cite{Yb}, and that we should be able  to develop some kind of understanding similar to \cite{ADY}. However, such a proof eludes us for now; it is not clear that a time-reversal argument may be the missing stone\dots \,Nevertheless, we present some further extensions of Bougerol's identity different from the ones found in the volume \cite{Yd}, which, hopefully should lead us in the future to a better understanding of Theorem \ref{theo1}.

\section{Discussion of and some theorems closely related to Theorem 1}\label{discussion}

\noindent{\bf (2.1)} Firstly we note that we may rewrite the identity in law \eqref{ids} in a seemingly slightly weaker form
\begin{align}\label{abs}
  (\sinh(|B_t|)\,,\,\sinh(L_t))\, \law \ (|\beta|(A_t)\,,\,\exp(-B_t)\lambda(A_t))\, \law \ (\exp(-B_t)|\beta|(A_t)\,,\,\lambda(A_t))\,.
  \end{align}
  This induces no loss of generality,
since the left-hand sides (without absolute values) of the expressions in \eqref{ids} only differ from the ones with absolute values in \eqref{abs} by multiplication by a symmetric Bernoulli variable, independent of the remaining quantities.

Secondly, it is well-known that the law of the two-dimensional vector $(|B_t|,L_t)$ is symmetric. More precisely,  it is given by:
\begin{align}\label{btlt}
\P(|B_t|\in dx\,,\,L_t\in d\ell)=\frac{2 (x+\ell)}{\sqrt{2\pi
    t^3}}\exp\lf-\frac{(x+\ell)^2}{2t}\ri dx\,d\ell\,,\qquad x\geq 0, \ell\geq0\,,
\end{align}
as it can be checked from L\'evy's identity (Theorem VI.2.3 in  \cite{RY}, p.240)
and the reflexion principle (Exercise 3.14 in  \cite{RY}, p.110, or Proposition 2.8.1 in \cite{KS}, p. 95). Hence, the common law of \eqref{abs} 
is also symmetric.

\noindent{\bf(2.2)} We now consider the right-hand sides of the first and third vectors in \eqref{ids} (or \eqref{abs}), and we deduce therefrom
\begin{align}\label{sinh}
\sinh(L_t)\ \law\ \lambda(A_t)\,.
\end{align}
Now, unlike for Bougerol's identity \eqref{1}, for which the possibility of an identity in law between processes is immediatley ruled out, since the left-hand side of \eqref{1} is not a martingale, whereas the right-hand side is, when one considers \eqref{sinh} it seems reasonable to wonder whether this identity might be valid at the level of the two increasing processes involved. However, the recent results in \cite{BY} also rule out this possibility. In fact, it is this uncertainty which prevented us from publishing an earlier version of this paper, as in \cite{DY}.

We point out also that for each fixed $t\geq 0$, we deduce from \eqref{eqproc} the rather puzzling identity in law
$$\lambda(A_t)\law  \sinh \left( \int_0^t\exp(-B_s) d\lambda(A_s)\right)\,,$$
which complements \eqref{sinh}. 

 \noindent{\bf(2.3)} {\em Reformulation in terms of Brownian suprema}.
A celebrated identity in distribution due to Paul L\'evy states that
$$((\overline B_t-B_t,\overline B_t)\,,t\ge0) \ \law \ ((|B_t|\,,L_t)\,,t\ge0)$$ 
where $\overline B_t=\sup_{s\le t}B_s$ denotes the supremum of the Brownian trajectory up to time $t$. This enables us to reformulate \eqref{abs} in the form
\begin{eqnarray*}
  (\sinh(\overline B_t-B_t)\,,\,\sinh(\overline B_t))\, &\law& \ ((\overline{\beta}-\beta)(A_t)\,,\,\exp(-B_t)\overline{\beta}(A_t))\\
  &\law& \ (\exp(-B_t)(\overline{\beta}-\beta)(A_t)\,,\,\overline{\beta}(A_t))\,,
\end{eqnarray*}
with $\overline \beta(t)=\sup_{s\le t}\beta(s)$. We leave to the interested reader further alternative reformulations of this identity in the same vein.

\noindent{\bf(2.4)} {\em A partial interpretation in terms of the Bessel clock}.

We now discuss Bougerol's identity \eqref{1} in terms of a two-dimensional Bessel process. 
Specifically, let $(R_h,h\geq0)$ denote  2-dimensional Bessel process 
starting from 1, and let
 $$
 H_u=\dis\int_0^u\frac{dh}{R^2_h}\,, \qquad u\geq0\,,
 $$  
 the clock associated with $R$.
 The well-known skew-product decomposition of planar Brownian motion (\cite{IM}, p.270; \cite{M1};\cite{M2}) shows that the clock $H$ can be viewed as the inverse of the exponential Brownian functional $A$; consequently, considering the inverses of the increasing processes involved in
\eqref{sinh}, we obtain the following:

\begin{corollary}\label{C1}
Let $(\sigma_t,t\geq0)$ denote a
stable $(1/2)$ subordinator, precisely:
\begin{equation*}
    \sigma_t:=\inf\{u:\lambda_u\geq t\}\,,\;t\geq0\,,
\end{equation*}
independent from $(R_h,h\geq0)$. Then, for fixed $s$, one has:
\begin{equation}\label{7}
    H_{\sigma_s}\law\sigma_{a(s)}
\end{equation}
where $a(s)\equiv\arg \sinh (s)$.
\end{corollary}

It is interesting to compare Corollary \ref{C1} with the following consequence of \eqref{eqproc}:
\begin{corollary}\label{C2}There is the identity in law between processes
$$\left(\sigma_t\right)_{t\geq 0}\law \left( H_{\sigma_{\eta(t)}}\right)_{t\geq 0}
$$
where $\sigma$ is as in Corollary \ref{C1} and $\eta: [0,\infty)\to [0,\infty)$ is the inverse bijection of the continuous strictly increasing process 
$$s\mapsto \int_0^s \frac{du}{R_{\sigma_u}}\,.$$
\end{corollary}
\proof Indeed, \eqref{eqproc} yields in terms of the Bessel clock 
$$\left( L_t
\right)_{t\geq 0} \law \left( \int_0^{A_t}\frac{d \lambda(u)}{R_u}\right)_{t\geq 0}\,.
$$
Our statement now follows from the easy fact that the process $\left(\sigma_{\eta(t)}, t\geq 0\right)$ is the right-inverse of the continuous increasing process $s\mapsto \int_0^s \frac{d\lambda(u)}{R_{u}}$. \QED

Again, one may wonder whether  identity (\ref{7}) holds at the level of increasing processes, however
the results in  \cite{BY} rule out this possibility. It is then natural to ask for which functionals
$\Phi: {\mathcal C}^{\uparrow}\to \R$ the identity
\begin{equation}\label{z}
\E(\Phi(H_{\sigma_{\cdot}}))= \E(\Phi(\sigma_{a(\cdot)}))
\end{equation}
may hold, where ${\mathcal C}^{\uparrow}$ stands for the space of c\`adl\`ag increasing paths $\omega: \R_+\to \R_+$. 
An element of positive response is provided by the following result.

\begin{theorem}\label{theo2}
Consider a measurable function $\Gamma: \R_+^3\to \R_+$ with $\Gamma(\cdot, 0, \cdot)=0$ and define
$$\Phi(\omega)= \sum_{s\geq 0} \Gamma(\omega_{s-}, \Delta \omega_s, s)\,,\qquad \omega\in {\mathcal C}^{\uparrow}\,,$$
where $\Delta \omega_s=\omega_s-\omega_{s-}$. Then \eqref{z} holds.

More precisely, if $\Gamma(x,y, s)= f(x,  s) g(y)$ for some measurable nonnegative functions $f$ and $g$ with $g(0)=0$, then 
$$
\E(\Phi(H_{\sigma_{\cdot}}))= \E(\Phi(\sigma_{a(\cdot)}))= C(f) D(g)\,,
$$
with
$$
C(f) = \int_0^{\infty} {d \lambda\over\sqrt{1+\lambda^2}}\E[f(H_{\sigma_\lambda},\lambda)]\ =\ 
 \int_0^{\infty} {d \lambda\over\sqrt{1+\lambda^2}}\E[f(\sigma_{a(\lambda)},\lambda)]
$$
and 
$$D(g) = \int_0^\infty\frac{dt}{\sqrt{2\pi t^3}}g(t)\,. $$
\end{theorem}

We observe that  Corollary \ref{C2} implies that the range of subordinated clock $H_{\sigma}$
has the same distribution as the range of subordinator  $\sigma$ (and hence is a regenerative set).
In particular we see that \eqref{z} holds whenever for a generic increasing path $\omega$,  $\Phi(\omega)$ only depends on the range of $\omega$. This provides a quick check of the identity 
$\E(\Phi(H_{\sigma_{\cdot}}))= \E(\Phi(\sigma_{a(\cdot)}))$
 in the special case when the function $\Gamma$ does not depend on the time parameter, i.e. $\Gamma(x,y,s)= \Gamma(x,y)$. 

\noindent{\bf(2.5)} {\em An amplification and a variant of Theorem 2}.

\begin{theorem}\label{theo3}
Let $a,b\geq0$ and  $\Gamma: \R_+^2\to \R_+$  a measurable function with $\Gamma(\cdot, 0)=0$. Introduce
\begin{equation*}
    \mathcal{H}^{a,b}(\Gamma)(\ell)=\E\left[\sum_{\lambda\leq\ell}(R_{\sigma_{\lambda-}})^a\frac{\Gamma(H_{\sigma_{\lambda-}},H_{\sigma_{\lambda}}-H_{\sigma_{\lambda-}})}{R^b_{\sigma_\lambda}}\right].
\end{equation*}

The following formula holds for $\Gamma =f\otimes g$:

\begin{equation*}
    \mathcal{H}^{a,b}(f\otimes g)(\ell)=h^-_{a-b}(f,\ell)h^+_b(g)\;,
\end{equation*}

where:
\begin{eqnarray*}
  h^-_c(f,\ell) &=& \int_0^\infty dt\,f(t)\E\left[\frac{e^{(c+1)B_t}}{\sqrt{2\pi A_t}}(1-e^{-\ell^2/2A_t})\right] \\
  h^+_b(g) &=& \int_0^\infty\frac{dt\,g(t)}{\sqrt{2\pi}}\E\left[\frac{e^{(2-b)B_t}}{A_t^{3/2}}\right].
\end{eqnarray*}

\end{theorem}

It may be interesting to point out that these formulas become simpler in the special case when $a=0$ and $b=1$. Indeed, one gets using \eqref{4} that
$$ h^+_1(g)=  \int_0^\infty\frac{dt}{\sqrt{2\pi t^3}}\, g(t)$$
and  then, using \eqref{3} and \eqref{2} that
$$  h^-_{-1}(f,\ell) =  \int_0^\infty\frac{dt\,f(t)}{\sqrt{2\pi t}}\left( 1- a'(\ell) \exp\left(-\frac{a^2(\ell)}{2t}\right) \right)$$
where $a$ and $a'$ are defined below \eqref{2}. 

We also stress that the quantities 
\begin{equation} \label{mpq}
m_{p,q}(t)=\E\left[ \frac{\exp(pB_t)}{A_t^q}\right] 
\end{equation}
arising in Theorem \ref{theo3}, have been studied in \cite{CMY} and \cite{Dufb, Dufc}.

\noindent{\bf(2.6)} {\em A variant of Theorem 2 which involves the windings of planar Brownian motion}.

The following variant of Theorem \ref{theo2} bears upon
a relationship between the continuous winding process of planar
Brownian motion, subordinated with $(\sigma_\lambda,
\lambda\geq0)$, and the standard Cauchy process.

\begin{theorem}\label{theo4}

Let $Z_u=|Z_u|\exp(i\theta_u)$, $u\geq0$, denote complex valued
Brownian motion, starting from: $1+i0$, with $(\theta_u,u\geq0)$
its continuous winding process. 
Let $(\sigma_\lambda,\lambda\geq0)$ denote the inverse local time process of a linear Brownian motion, so $\sigma$ is a stable(1/2) subordinator, which is assumed
to be independent from $(Z_u,u\geq0)$.
Finally, let  $(C_\alpha,\alpha\geq0)$ be a standard Cauchy process.\\

For any measurable $\Gamma:\R\times\R\to\R_+$ with $\Gamma(x,0)=0$, we have
$$
  \E\left[\sum_{\lambda\leq\ell} \Gamma(\theta_{\sigma_{\lambda-}},\theta_{\sigma_\lambda}-\theta_{\sigma_\lambda-})\right] =
 \E\left[\sum_{\lambda\leq\ell} \Gamma(C_{a(\lambda)-},C_{a(\lambda)}-C_{a(\lambda)-})\right]
$$
for all $\ell$'s. In particular, for fixed $\ell\geq0$, there is the equality in
law:

\begin{equation}\label{eqNEUF}
    \theta_{\sigma_\ell}\law C_{a(\ell)}\,.
\end{equation}

\end{theorem}

The reader interested in some applications of these
identities in law to functionals of the winding process
$(\theta_u,u\geq0)$ may refer to Vakeroudis \cite{V}. In
particular, the identity \eqref{eqNEUF} allows to apply D. Williams'
pinching method to yield yet another proof of Spitzer's celebrated
theorem:

\begin{equation}\label{eqDIX}
   \frac{2}{\log t}\;\theta_t\underset{t\to\infty}{\overset{(\textrm{law})}{\longrightarrow}}C_1\,.
\end{equation}

\noindent{\bf(2.7)} {\em The joint Laplace-Mellin transform of $(H_{\sigma_\lambda}\,,R_{\sigma_\lambda})$}.

We now come back to Theorem \ref{theo2}, or rather we discuss part of its proof, as given in paragraph (3.2) below. A by product of Lemma 1 therein is:
\begin{equation}\label{9}
    \E\lff\frac{1}{R_{\sigma_\lambda}}\,\bigg\vert\, H_{\sigma_\lambda}=h\rii=\frac{1}{\sqrt{1+\lambda^2}}\,,\end{equation}
an intriguing identity, which made us suspect for a moment that
$R_{\sigma_\lambda}$ and $H_{\sigma_\lambda}$ might be
independent. This is not the case, as we discovered by computing
the joint Laplace-Mellin transform of
$(H_{\sigma_\lambda},R_{\sigma_\lambda})$:

\begin{theorem}\label{theo4}
The following formulae hold:

\begin{eqnarray}
  \E\lff\frac{1}{(R_{\sigma_\lambda})^{2b}}\exp(-\frac{\mu^2}{2}H_{\sigma_\lambda})\rii&=& \E^{(\mu)}\lff\frac{1}{(R_{\sigma_\lambda})^{2b+\mu}}\rii \nonumber\\
  &=& C_{b,\mu}Ê\times \frac{F\lf\frac{\mu+1}{2}-b,\frac{\mu}{2}+1-b,\mu+1;-1/{\lambda^2}\ri}{(1+\lambda^2)^{2b-1/2}(\lambda^2)^{\frac{\mu+1}{2}-b}}\,,\qquad \label{10}
\end{eqnarray}
where $F\equiv {}_2F_1$ denotes the classical hypergeometric family of functions with three parameters, and $\E^{(\mu)}$ refers to the expectation with respect to  the probability measure
$\P^{(\mu)}$ under which $(R_t, t\geq0)$ is a Bessel process with index $\mu$ (i.e. of dimension $2+2\mu$) and  started from $R_0=1$,
and
$$C_{b,\mu} = \frac{\Gamma(b+\frac{\mu}{2}+\frac{1}{2}) \Gamma(1+\frac{\mu}{2}-b)}
{\Gamma(\frac{1}{2}) \Gamma(1+\mu)}.$$\end{theorem}
As a partial check for formula \eqref{10}, we have made verifications with $b=1/2$, $b=-\mu/2$ (the result should be $1$), and $b=0$ (the result is $\exp(-\mu a(\lambda))=(\lambda+\sqrt{1+\lambda^2})^{-\mu}$).
Let us give some details for $b=1/2$:
We note that, for $b=1/2$, (\ref{10}) simplifies, as in the numerator
\begin{equation*}
    F\lf\frac{\mu}{2},\frac{\mu}{2}+\frac{1}{2},\mu+1;-y\ri=\lf\frac{2}{1+\sqrt{1+y}}\ri^\mu\;,
\end{equation*}
and in the denominator
$$
(1+\lambda^2)^{1/2}(\lambda^2)^{\mu/2}\ \equiv\ \lambda^\mu(1+\lambda^2)^{1/2}.
$$
Hence, using the fact that $C_{1/2, \mu}= 2^{-\mu}$, formula \eqref{10} simplifies to:
$$2^{-\mu} \cdot
{2^{\mu}\over\left(1+\sqrt{1+{1\over\lambda^2}}\right)^\mu}\cdot {1\over \lambda^\mu\sqrt{1+\lambda^2}}  \equiv\ {1\over (\sqrt{1+\lambda^2}+\lambda)^\mu\sqrt{1+\lambda^2}}\,,
$$
which confirms identity \eqref{9}, since the previous RHS expression equals 
$$\frac{1}{\sqrt{1+\lambda^2}}\E\left(\exp\left(-\frac{\mu^2}{2} H_{\sigma_{\lambda}}\right)Ê\right)\,.$$

\section{Proof of  Theorem 1}\label{proof1}

\textbf{(3.1)} We start by recalling some well-known facts about Brownian motion running up to an independent exponential time, which will be useful for the proof. 
In this subsection,  $S_p$ denotes an
exponential random variable with parameter
$p>0$, independent from the Brownian motion $B$.
For $\tg$, denote $g_t=\sup\{u<t:B_u=0\}$ the last zero of $B$
before $t$. It is known that the processes $(B_u,u\leq g_{S_p})$ and $(B_{g_{S_p}+u},u\leq S_p-g_{S_p})$ are independent. As a consequence, the variables $L_{S_p}(\equiv L_{g_{S_p}})$ and
$B_{S_p}$ are independent. Moreover, since $L_t$ and $|B_t|$ have the same law (see \eqref{btlt}), the same applies to $L_{S_p}$ and $|B_{S_p}|$. Their common density is
\begin{equation*}
    \sqrt{2p}\exp(-\sqrt{2p}u),\qquad u\ge0
\end{equation*}
(this is because $\P(L_{S_p}\geq\ell)=\P(S_p\geq\tau_\ell)=\E[\exp(-p\tau_\ell)]=\exp(-\ell\sqrt{2p})$, if $\tau_\ell$ is the time $L_\cdot$ reaches $\ell$).
 An equivalent way to express this property is
$$
\sqrt{2\mathbf{e}}(|\beta(1)|,\lambda(1))\ \law\ (\mathbf{e},\mathbf{e}'),
$$
where on the left $\mathbf{e}\law S_1$ is independent of $\beta$, and on the right the two variables are independent copies of $S_1$.

\noindent\textbf{(3.2)} Recall the discussion in paragraph (2.1). Our main goal is to show that 
\begin{equation}\label{eq25}
    (\sinh(|B_t|),\sinh(L_t))\law(\exp(-B_t)\sqrt{A_t}|\beta(1)|,\sqrt{A_t}\lambda(1))\,.
\end{equation}
This will be done by computing the joint Mellin transforms on either side, but before doing so we replace $t$ with an exponential time $S_p$ and multiply both sides by $\sqrt{2\mathbf{e}}$, assuming implicitly that $S_p$, $\mathbf{e}$,  $B$ and $\beta$ are independent. What will be proved is:
\begin{equation}\label{eq25a}
   \sqrt{2\mathbf{e}} (\sinh(|B_{S_p}|),\sinh(L_{S_p}))\law\sqrt{2\mathbf{e}} (\exp(-B_{S_p})\sqrt{A_{S_p}}|\beta(1)|,\sqrt{A_{S_p}}\lambda(1))\,.
\end{equation}

From the one-dimensional Bougerol identity \eqref{1}, the left-hand side of \eqref{eq25a} has the same distribution as
\begin{equation}\label{e1}
    \sqrt{2\mathbf{e}}(|N|\sqrt{A_{S_p}},|N'|\sqrt{A'_{S'_p}})\,,
\end{equation}
where $N$, $N'$ are independent standard normals and $A'$ a copy of $A$ which is also independent of the other quantities. On the right-hand side of \eqref{eq25a}, use the facts in paragraph (3.1) to obtain
\begin{align}\label{e2}
\sqrt{2\mathbf{e}} (\exp(-B_{S_p})\sqrt{A_{S_p}}|\beta(1)|,\sqrt{A_{S_p}}\lambda(1))&\law
(\exp(-B_{S_p})\sqrt{A_{S_p}}\mathbf{e}\,,\sqrt{A_{S_p}}\mathbf{e}')\nonumber\\
&\law (\exp(-B_{S_p})\sqrt{A_{S_p}}|N|\sqrt{2\mathbf{e}}\,,\sqrt{A_{S_p}}|N'|\sqrt{2\mathbf{e}'})\,.
\end{align}
Squaring, we are left with calculating the joint Mellin transforms of 
\begin{align}\label{e3}
\mathbf{e}(A_{S_p},{A'_{S'_p}})\quad\text{and}\quad (\exp(-2B_{S_p}){A_{S_p}}{\mathbf{e}}\,,{A_{S_p}}{\mathbf{e}'})\,,
\end{align}
and verifying that they are equal (the $N,N'$ on both sides of \eqref{e1} and \eqref{e2} can be cancelled). 

The essential ingredient for these Mellin transforms is (see \cite{Ya}, paper$\,\sharp$6,
p. 94)
\begin{equation}\label{eq12}
    A^{(\nu)}_{S_p}\law\frac{\beta_{1,a}}{2\gamma_b}\;,
\end{equation}
where $A^{(\nu)}$ denotes the exponential functional
$$   A^{(\nu)}_t=\int_0^tds\exp(2B^{(\nu)}_s)\,,\qquad t\ge0,
$$
of a Brownian motion with drift $\nu$, $B_s^{(\nu)}= B_s + \nu s$,
 $\beta_{u,v}$ is a beta variable with parameters $(u,v)$,
 $\gamma_b$ is a gamma variable with parameter $b$, and
 $$
  a= a(\nu,p)=\frac{1}{2}(\nu+\sqrt{2p+\nu^2})\,,\qquad b=b(\nu,p)=\frac{1}{2}(-\nu+\sqrt{2p+\nu^2})\,.
  $$The Mellin transform of $A^{(\nu)}_{S_p}$ is then
  $$
  \E(A^{(\nu)}_{S_p})^r= 2^{-r}{\Gamma(1+a)\Gamma(1+r)\Gamma(b-r)\over \Gamma(1+a+r)\Gamma(b)}\,.
  $$
On the one hand, since $a(0,p)=b(0,p)=\sqrt{p\over2}$,
\begin{align}\label{other1}
\E[(\mathbf{e}A_{S_p})^c(\mathbf{e}A'_{S'_p})^d]&=2^{-c-d} {\Gamma(1+c+d)\Gamma(1+a)\Gamma(1+c)\Gamma(b-c)\Gamma(1+a)\Gamma(1+d)\Gamma(b-d) \over \Gamma(1+a+c)\Gamma(1+a+d)\Gamma(b)^2}\nonumber\\
&= 2^{-c-d}\,{p\over2}\  {\Gamma(1+c+d)\Gamma(1+c)\Gamma(\sqrt{p\over2}-c)\Gamma(1+d)\Gamma(\sqrt{p\over2}-d) \over \Gamma(1 +c+ \sqrt{p\over2})\Gamma(1 +d+ \sqrt{p\over2})}\,.
\end{align}
On the other hand, recalling the definition of $A^{(\nu)}_t$ and using the Girsanov-Cameron-Martin theorem, we obtain
\begin{align*}
\E[(\exp(-2B_{S_p}){A_{S_p}}{\mathbf{e}})^c({A_{S_p}}{\mathbf{e}'})^d]&= 
\Gamma(1+c)\Gamma(1+d)\E [\exp(-2cB_{S_p})(A_{S_p})^{c+d}]\\
&=\Gamma(1+c)\Gamma(1+d)\E [\exp(2c^2{S_p})(A^{(-2c)}_{S_p})^{c+d}].
\end{align*}
From the elementary identity
$$
\P(S_p\in ds\,;\, \exp(\eta S_p))={p\over p-\eta}\P(S_{p-\eta}\in ds)\,,\qquad \eta<p\,,
$$
we find, letting $q=p-2c^2$,
\begin{align}\label{other2}
&\E[(\exp(-2B_{S_p}){A_{S_p}}{\mathbf{e}})^c({A_{S_p}}{\mathbf{e}'})^d]\nonumber\\
&= 
{p\over q}\Gamma(1+c)\Gamma(1+d)\E [(A^{(-2c)}_{S_q})^{c+d}]\nonumber\\
&=2^{-c-d} {p\over q}\Gamma(1+c)\Gamma(1+d)
  {\Gamma(1+a(-2c,q))\Gamma(1+c+d)\Gamma(b(-2c,q)-c-d)\over \Gamma(1+a(-2c,q)+c+d)\Gamma(b(-2c,q))}\,.
\end{align}
Now, $a(-2c,q)=-c+\sqrt{p\over2}$, $b(-2c,q)=c+\sqrt{p\over2}$, and thus, comparing \eqref{other1}-\eqref{other2},
\begin{align*}
& {p\over q}{\Gamma(1+a(-2c,q))\Gamma(b(-2c,q)-c-d)\over \Gamma(1+a(-2c,q)+c+d)\Gamma(b(-2c,q))}\\
&=\ {p\over p-2c^2}{\Gamma(1-c+\sqrt{p\over2})\Gamma(-d+\sqrt{p\over2})\over
\Gamma(1+d+\sqrt{p\over2})\Gamma(c+\sqrt{p\over2})}\\
&= {p\over p-2c^2}{\textstyle(\sqrt{p\over2}-c)(\sqrt{p\over2}+c)}{\Gamma(\sqrt{p\over2}-c)\Gamma(\sqrt{p\over2}-d) \over \Gamma(1 +c+ \sqrt{p\over2})\Gamma(1 +d+ \sqrt{p\over2})}\\
&= {p\over 2}{\Gamma(\sqrt{p\over2}-c)\Gamma(\sqrt{p\over2}-d) \over \Gamma(1 +c+ \sqrt{p\over2})\Gamma(1 +d+ \sqrt{p\over2})}\,.
\end{align*}
Note that, along the way, it was necessary to assume $q=p-2c^2>0$, so that $c$ needed to be taken small enough, and likewise for $d$, precisely: $c,d<\sqrt{p\over2}$. But, even with these restrictions, we can conclude the proof of the identity in law of the two vectors in \eqref{e3}, thus ending the proof of \eqref{eq25}.

\noindent\textbf{(3.3)} The second identity in Theorem 1 may be proved rather simply, by first noting that
\begin{equation}\label{eq8}
    (\beta(A_t),\exp(-B_t)\lambda(A_t))\law(\sqrt{A_t}\beta(1),\exp(-B_t)\sqrt{A_t}\lambda(1))
\end{equation}
and then recalling (from the proof of Bougerol's identity in \cite{ADY}) that time reversal 
\begin{equation*}
    (B_t-B_{(t-u)}\,,\,0\le u\leq t)\law(B_u\,,\,0\le u\leq t)
\end{equation*}
implies
\begin{align*}
(A_t\,,\,\exp(-2B_t)A_t)\law(\exp(-2B_t)A_t\,,\,A_t)\,.
\end{align*}
This completes the proof of Theorem \ref{theo1}. 

\section{Proofs of  Theorems 2 to 5}\label{proof6}
\noindent\textbf{(4.1)} {\em Proof of Theorem 2}.
 (a) A key for the proof of Theorem 2 is the following interesting, and puzzling, identity, as discussed in Subsection 2.7.

\begin{lemma} For any measurable function $f:\R_+\mapsto \R_+$, and any $s\ge0$, we have:
$$
\E\left[{1\over R_{\sigma_s}}f(H_{\sigma_s})\right]={1\over\sqrt{1+s^2}}\E[f(\sigma_{a(s)})].
$$
\end{lemma}

\proof
For all $q,t, \varepsilon >0$, we have from \eqref{7} that
$$\E(\exp (-q \sigma_{a(t)})-\exp(-q \sigma_{a(t+\varepsilon)})) = 
\E(\exp (-q H_{\sigma_t})-\exp(-q H_{\sigma_{t+\varepsilon}}))\,.$$
The LHS can be computed explicitly and we obtain
$$\exp(-a(t)\sqrt{2q})-\exp(-a(t+\varepsilon)\sqrt{2q})) \sim \frac{\varepsilon}{\sqrt{1+t^2}} \sqrt{2q}Ê\exp(-a(t)\sqrt{2q})\,, \quad \varepsilon \to 0\,.$$

We next turn our attention to the RHS and apply the Markov property. In this direction, it is convenient to introduce a two-dimensional Bessel process $R'$ which is independent of $R$ and write $H'$ for its clock. 
Likewise, $\sigma'$ refers to an independent subordinator which has the same distribution as $\sigma$. 
For every $r>0$, the notation $\P'_r$ refers to 
the law under which $R'_0=r$ and $\E'_r$ to the mathematical expectation under $\P'_r$. We point out that the scaling property implies the identities
$$\E'_r\left(1-\exp (-q H'_{\sigma'_{\varepsilon}})\right)=
\E'_1\left(1-\exp (-q H'_{r^{-2}Ñ\sigma'_{\varepsilon}})\right)
=\E'_1\left(1-\exp (-q H'_{\sigma'_{\varepsilon/r}})\right).
$$
Of course we can also express the RHS as $\E\left(1-\exp (-q H_{\sigma_{\varepsilon/r}})\right)$. 
Writing $\sigma_{t+\varepsilon}= \sigma_t + \sigma'_{\varepsilon}$, with $\sigma'_{\varepsilon}$ independent from $\sigma_t$ and $R$, we get from an application of the Markov property
\begin{eqnarray*}
\E(\exp (-q H_{\sigma_t})-\exp(-q H_{\sigma_{t+\varepsilon}}))
&=& \E\left(\exp(-q H_{\sigma_t}) \E'_{R_{\sigma_t}}\left(1-\exp (-q H'_{\sigma'_{\varepsilon}})\right) \right) \\
&=& \E\left (\exp(-q H_{\sigma_t}) \E'_{1}\left(1-\exp (-q H'_{\sigma'_{\varepsilon/R_{\sigma_t}}})\right) \right)\\
&=& \E\left(\exp(-q H_{\sigma_t}) \left(1-\exp\{-a(\varepsilon/R_{\sigma_t})\sqrt{2q}\} \right) \right)
\end{eqnarray*}
where the third equality stems from \eqref{7}. Note that when $\varepsilon \to 0+$, the preceding quantity is equivalent to
$$\varepsilon  \sqrt{2q} \E\left(\exp(-q H_{\sigma_t}) \frac{1}{R_{\sigma_t}}\right)\,.$$

Putting the pieces together, we arrive at 
$$\E\left (\exp(-q H_{\sigma_t}) \frac{1}{R_{\sigma_t}}\right)=\frac{1}{\sqrt{1+t^2}} Ê\exp(-a(t)\sqrt{2q}) =\frac{1}{\sqrt{1+t^2}} \E(\exp (-q \sigma_{a(t)}))$$
for all $t,q>0$, which establishes Lemma 1.\QED 

\noindent (b) To end the proof of Theorem 2, we introduce the following notation concerning jump intensity measures:
\begin{eqnarray}
  && \mathcal{H}(\Gamma)(s)=\E\left [\sum_{\lambda\leq s}\Gamma (H_{\sigma_{\lambda^-}},H_{\sigma_{\lambda}}-H_{\sigma_{\lambda^-}})\right ] \nonumber\\
  \textrm{and} && \label{8}\\
  && \mathcal{K}(\Gamma )(s)=\E\left [\sum_{\alpha\leq a(s)}\Gamma (\sigma_{\alpha^-},\sigma_{\alpha}-\sigma_{\alpha^-})\right ] \nonumber
\end{eqnarray}
for given $s$, and $\Gamma :\RR_+\times\RR_+\to\RR_+$, Borel, such that $\Gamma(x,0)=0$.
Then, in order to finish the proof of Theorem 2, it suffices to take: $\Gamma=f\otimes g$ and to show the following:
\begin{eqnarray*}
  \mathcal{H}(f\otimes g)(s) &=& h(f)(s)\int_0^\infty\frac{dt}{\sqrt{2\pi t^3}}g(t) \\
  \mathcal{K}(f\otimes g)(s) &=& k(f)(s)\int_0^\infty\frac{dt}{\sqrt{2\pi t^3}}g(t)
\end{eqnarray*}
with, furthermore, the quantities $h(f)(s)$ and $k(f)(s)$ being equal, and equal to:
\begin{eqnarray*}
  h(f)(s) &=& \int_0^sd\lambda \E \lff\frac{1}{R_{\sigma_\lambda}}f(H_{\sigma_\lambda})\rii=\int_0^s\frac{d\lambda}{\sqrt{1+\lambda^2}}\E[f(H_{\sigma_\lambda})] \\
  \parallel && \\
  k(f)(s) &=& \int_0^{a(s)}d\alpha \E[f(\sigma_\alpha)]=\int_0^s\frac{d\lambda}{\sqrt{1+\lambda^2}}\E[f(\sigma_{a(\lambda)})]\,.
\end{eqnarray*}

Now, concerning $\mathcal{K}(\Gamma)(s)$,  since the Lévy
measure of the subordinator $(\sigma_\alpha,\alpha\geq0)$ is
$\dis\frac{dt}{\sqrt{2\pi t^3}}$, we have:
\begin{eqnarray*}
  \mathcal{K}(\Gamma)(s) &=& \E\lff\int_0^{a(s)}d\alpha f(\sigma_{\alpha-})\rii\int_0^\infty\frac{dt}{\sqrt{2\pi t^3}}g(t) \\
  &=& \int_0^s\frac{du}{\sqrt{1+u^2}}\E[f(\sigma_{a(u)})]\int_0^\infty\frac{dt}{\sqrt{2\pi
  t^3}}g(t)\,.
\end{eqnarray*}

Concerning $\mathcal{H}(\Gamma)(s)$, starting again with the same
argument (i.e: the knowledge of the L\'evy measure of
$(\sigma_u,u\geq0)$), we obtain:
\begin{eqnarray}
  \mathcal{H}(\Gamma)(s) &=& \E\lff\int_0^sd\lambda\int_0^\infty\frac{dt}{\sqrt{2\pi t^3}}f(H_{\sigma_\lambda})g(H_{\sigma_\lambda+t}-H_{\sigma_\lambda})\rii\nonumber \\
  &=& \E\lff\int_0^sd\lambda\int_0^\infty\frac{dt}{\sqrt{2\pi t^3}}f(H_{\sigma_\lambda})\E'_{R_{\sigma_\lambda}}(g(H'_t))\rii\label{a}\
\end{eqnarray}
(from the Markov property for $R$).
However, by scaling, we have:
\begin{equation}\label{b}
    \E'_\rho[g(H'_t)]=\E[g(H_{t/\rho^2})]
\end{equation}
so that, plugging (\ref{b}) in (\ref{a}), we obtain:
\begin{eqnarray}
  \int_0^\infty\frac{dt}{\sqrt{2\pi t^3}}\E'_\rho[g(H'_t)] &=& \int_0^\infty\frac{dt}{\sqrt{2\pi t^3}}\E[g(H_{t/\rho^2})]\nonumber \\
  &=& \frac{1}{\rho}\int_0^\infty\frac{du}{\sqrt{2\pi u^3}}\E[g(H_u)]\,.\label{c}
\end{eqnarray}
Note that, since the inverse of $\{u\to H_u\}$ is: $t\to
A_t=\dis\int_0^tdve^{2B_v}$, we have:
\begin{eqnarray*}
  \int_0^\infty\frac{du}{\sqrt{2\pi u^3}}\E[g(H_u)] &=& \frac{1}{\sqrt{2\pi}}\E\lff\int_0^\infty dt g(t)\frac{e^{2B_t}}{\sqrt{A_t^3}}\rii \ =\  \int_0^\infty dt g(t)\frac{1}{\sqrt{2\pi t^3}}
\end{eqnarray*}
by (\ref{4}).

Going back to (\ref{a}), we have obtained:
\begin{equation}\label{d}
  \mathcal{H}(\Gamma)(s)=\E\lff\int_0^sd\lambda\frac{1}{R_{\sigma_\lambda}}f(H_{\sigma_\lambda})\rii\int_0^\infty\frac{dt}{\sqrt{2\pi t^3}}g(t)
\end{equation}
b) Finally, to obtain the equality between $\mathcal{H}(\Gamma)(s)$ and
$\mathcal{K}(\Gamma)(s)$, it remains to show, with the notation in the
statement of Theorem \ref{theo3}, that:
\begin{equation*}
    h(f)(s)=k(f)(s)\,,\;\textrm{for every $f\geq0$, Borel.}
\end{equation*}
Again, it suffices to prove this for $f_\lambda(a)=e^{-\lambda
a}$, for any $\lambda\geq0$.
Now we have:
\begin{eqnarray*}
  \E[\exp(-\nu H_{\sigma_s})] &=& 1+\E\left[\sum_{\theta\leq s}\lf e^{-\nu H_{\sigma_\lambda}}-e^{-\nu H_{\sigma_{\lambda^-}}}\ri\right] \\
  &=& 1+\E\left[\sum_{\theta\leq s}e^{-\nu H_{\sigma_{\lambda^-}}}\lf e^{-\nu(H_{\sigma_\lambda}-H_{\sigma_{\lambda^-}})}-1\ri\right] \\
  &=& 1+\mathcal{H}(f_\nu\otimes g_\nu)(s) \\
  (\textrm{where}&:& f_\nu(a)=\exp(-\nu a);g_\nu(b)=(e^{-\nu b}-1)) \\
  &=& 1+h(f_\nu)(s)\lf\int_0^\infty\frac{dt}{\sqrt{2\pi t^3}}g_\nu(t)\ri.
\end{eqnarray*}
On the other hand:
\begin{equation*}
    \E[\exp(-\nu\sigma_{a(s)})]=1+\int_0^{a(s)}d\alpha
    f_\nu(\alpha)\int_0^\infty\frac{dt}{\sqrt{2\pi t^3}}g_\nu(t)\,.
\end{equation*}
Thus, explicitly:
\begin{equation}\label{46}
    \E[\exp(-\nu H_{\sigma_s})]=1-(h(f_\nu)(s))\sqrt{2\nu}
\end{equation}
whereas:
\begin{equation}\label{47}
    \E[\exp(-\nu\sigma_{a(s)})]=1-k(f_\nu)(s)\sqrt{2\nu}\,.
\end{equation}

Since the left hand sides of (\ref{46}) and (\ref{47}) are equal,
so are the right hand sides, therefore:

\begin{equation}\label{48}
    \forall f\geq0\,,\;\textrm{Borel}\,,\;h(f)(s)=k(f)(s)\,.
\end{equation}
Hence, in complete generality:
\begin{equation}\label{49}
    \mathcal{H}(\Gamma)(s)=\mathcal{K}(\Gamma)(s)\,,
\end{equation}
which finishes the proof of Theorem 2.

\noindent\textbf{(4.2)} {\em  Proof of Theorem 3}. Here are the main steps of this proof, which is quite similar to that of Theorem 2:

1) We first transform
\begin{eqnarray*}
    \mathcal{H}^{a,b}(\Gamma)(\ell)&=&\E\left[\sum_{\lambda\leq\ell}(R_{\sigma_{\lambda-}})^af(H_{\sigma_{\lambda-}})\frac{g(H_{\sigma_{\lambda}}-H_{\sigma_{\lambda-}})1_{(\sigma_\lambda>\sigma_\lambda-)}}{(R_{\sigma_\lambda})^b}\right]\\
    &=&\E\left[\int_0^\ell d\lambda(R_{\sigma_\lambda})^af(H_{\sigma_{\lambda}})\int_0^\infty\frac{dt}{\sqrt{2\pi t^3}}\frac{g(H_{\sigma_{\lambda}+t}-H_{\sigma_\lambda})}{(R_{\sigma_\lambda+t})^b}\right] \\
    &=&\E\left[\int_0^\ell d\lambda(R_{\sigma_\lambda})^af(H_{\sigma_{\lambda}})E_{R_{\sigma_\lambda}}\left[\int_0^\infty\frac{dt}{\sqrt{2\pi t^3}}\frac{g(H_t)}{(R_t)^b}\right]\right].
\end{eqnarray*}

We begin by studying
\begin{eqnarray*}
    h^{(+)}(r,g)&=&\E_r\left[\int_0^\infty\frac{dt}{\sqrt{2\pi t^3}}\frac{g(H_t)}{R^b_t}\right]\\
    &=&\E\left[\int_0^\infty\frac{dt}{\sqrt{2\pi t^3}}\frac{g(H_{t/r^2})}{(rR_{t/r^2})^b}\right]\;\textrm{(by scaling)}\\
    &=&\frac{1}{r^b}\E\left[\int_0^\infty\frac{dt}{\sqrt{2\pi t^3}}\frac{g(H_{t/r^2})}{(R_{t/r^2})^b}\right]\\
    &=&\frac{1}{r^b}\lf\frac{1}{r}\ri \E\left[\int_0^\infty\frac{du}{\sqrt{2\pi u^3}}\frac{g(H_u)}{(R_u)^b}\right]\\
    &\equiv&\frac{1}{r^{b+1}}h^{(+)}_b(g)\,.
\end{eqnarray*}

We then study:
\begin{eqnarray*}
    h^{(+)}_b(g)&=&\E\left[\int_0^\infty\frac{du}{\sqrt{2\pi u^3}}\frac{g(H_u)}{(R_u)^b}\right]\\
    &=&\E\left[\int_0^\infty\frac{dte^{2B_t}}{\sqrt{2\pi A_t^3}}\frac{g(t)}{(\exp(bB_t))}\right]\\
    &=&\E\left[\int_0^\infty\frac{dte^{(2-b)B_t}}{\sqrt{2\pi A_t^3}}g(t)\right]\\
    &=&\,\int_0^\infty\frac{dtg(t)}{\sqrt{2\pi}}\E\left[\frac{e^{(2-b)B_t}}{\sqrt{A_t^3}}\right]\\
    &=&\int_0^\infty\frac{dtg(t)}{\sqrt{2\pi}}m_{2-b,3/2}(t)\,,
\end{eqnarray*}
where the quantity $m_{p,q}(t)$ has been defined in \eqref{mpq}. \\

3) Let us come back to
\begin{eqnarray*}
    \mathcal{H}^{a,b}(\Gamma)(\ell)&=&\E\left[\int_0^\ell d\lambda(R_{\sigma_{\lambda}})^{a-b-1}f(H_{\sigma_{\lambda}})\right]h^{(+)}_b(g)\\
    &=&h^{(-)}_{a-b}(f,\ell)h^{(+)}_b(g)\textrm{, for }F=f\otimes g\,.
\end{eqnarray*}
Thus, our next aim is to study:

\begin{equation*}
    h^{(-)}_c(f,\ell):=\E\left[\int_0^\ell d\lambda(R_{\sigma_{\lambda}})^{c-1}f(H_{\sigma_{\lambda}})\right].
\end{equation*}
We can re-express this quantity as 
\begin{eqnarray*}
   h^{(-)}_c(f,\ell)  &=& \E\left[\int_0^{\tau_\ell}dL_u(R_u)^{c-1}f(H_u)\right] \\
   &=&\int_0^\infty\frac{du}{\sqrt{2\pi u}}\P(L_u<\ell|B_u=0)\E[(R_u)^{c-1}f(H_u)] \\
  &=&\int_0^\infty\frac{du}{\sqrt{2\pi u}}\P(\sqrt{u}\sqrt{2\textbf{e}}<\ell)\E[(R_u)^{c-1}f(H_u)] \\
  &=&\int_0^\infty\frac{du}{\sqrt{2\pi u}}\P(\textbf{e}\leq\ell^2/2u)\E[(R_u)^{c-1}f(H_u)] \\
  &=&\int_0^\infty\frac{du}{\sqrt{2\pi u}}(1-\exp(-\ell^2/2u))\E[(R_u)^{c-1}f(H_u)] \\
  &=&\E\left[\int_0^\infty\frac{dA_t}{\sqrt{2\pi A_t}}(1-e^{-\ell^2/2A_t})\exp((c-1)B_t)f(t)\right] \\
  &=&\E\left[\int_0^\infty dt\frac{e^{(c+1)B_t}}{\sqrt{2\pi A_t}}(1-e^{-\ell^2/2A_t})f(t)\right] \\
  &=&\int_0^\infty dtf(t)\E\left[\frac{e^{(c+1)B_t}}{\sqrt{2\pi A_t}}(1-e^{-\ell^2/2A_t})\right] \,.
\end{eqnarray*}

\noindent\textbf{(4.3)} {\em Proof of Theorem 4}. It is a simple consequence of Theorem 2, once one uses the well-known skew product reresentation of $\theta_t=\gamma_{H_t}$, where $(\gamma_u,u\ge0)$ is a real-valued Brownian motion independent from $(H_t,t\ge0)$ (we already gave some references before Corollary 1). Then all one needs to do is to "freeze" $\gamma$ first, then apply Theorem 2, and finally use Spitzer's representation of the Cauchy process as 
$$
(C_\alpha,\alpha\ge0) \law (\gamma_{\sigma_\alpha},\alpha\ge0).
$$

\noindent\textbf{(4.4)} {\em Proof of Theorem 5}. The first equality follows from the (local) absolute continuity relationship between the laws of different Bessel processes, see, {\it e.g.}, \cite{Ya}.

Thus, it remains to prove the second equality. For this purpose, we use
the same arguments as in the proof of (16) in \cite{DY}; here are some details. 

Let
$$J\deff \E^{(\mu)}\lff\frac{1}{(R_{\sigma_\lambda})^{2b+\mu}}\rii
= \frac{1}{\Gamma(b+\frac{\mu}{2})} \int_0^{\infty} du \ u^{b+\frac{\mu}{2}-1} \E^{(\mu)}\lff\exp(-u R^2_{\sigma_\lambda})\rii\,.$$
There is a classical expression for $\E^{(\mu)}\lff\exp(-u R^2_t)\rii$ (see, e.g. \cite{RY} on page 441), which yields
$$\E^{(\mu)}\lff\exp(-u R^2_{\sigma_\lambda})\rii = \E\left(\left(1+2u\sigma_{\lambda}\right)^{-1-\mu}
\exp\left(-\frac{u}{1+2u\sigma_{\lambda}}\right)\right)\,.$$

Using Fubini and making the change of variables $v=2u\sigma_{\lambda}/(1+2u\sigma_{\lambda})$, we obtain the following expression for the integral
\begin{eqnarray*}I &\deff& \int_0^{\infty} du\ u^{b+\frac{\mu}{2}-1} \left(1+2u\sigma_{\lambda}\right)^{-1-\mu}
\exp\left(-\frac{u}{1+2u\sigma_{\lambda}}\right) \\
&=&  \frac{1}{(2\sigma_{\lambda})^{b+\mu/2}}\int_0^1dv \  v^{b+\frac{\mu}{2}-1} (1-v)^{\frac{\mu}{2}-b}Ê\exp\left( - \frac{v}{2\sigma_{\lambda}}\right)\,.
\end{eqnarray*}
Hence, using Fubini again, we obtain:
\begin{equation}
\label{last}
J=\frac{1}{\Gamma(b+\mu/2)}\int_0^1dv \  v^{b+\frac{\mu}{2}-1} (1-v)^{\frac{\mu}{2}-b}Ê
\E\left[ \frac{1}{(2\sigma_{\lambda})^{b+\mu/2}} \exp\left( - \frac{v}{2\sigma_{\lambda}}\right)\right]\,.
\end{equation}

To compute this last expectation, which we denote by $K$, we use 
$$ \frac{1}{2\sigma_{\lambda}}\law \frac{N^2}{2\lambda^2}\law \frac{\gamma_{1/2}}{\lambda^2}$$
where $\gamma_{1/2}$ is a standard gamma$(1/2$)-variable.
Then we obtain 
$$K= \frac{\Gamma(b+\frac{\mu}{2}+\frac{1}{2})}{\Gamma(\frac{1}{2})\lambda^{2b+\mu}(1+v\lambda^{-2})^{b+\frac{\mu}{2}+\frac{1}{2}}}\,.$$
Plugging this in \eqref{last}, we obtain: 
$$
J=\frac{\Gamma(b+\frac{\mu}{2}+\frac{1}{2})}{
\Gamma(\frac{1}{2})\Gamma(b+\mu/2)}\int_0^1dv \  v^{b+\frac{\mu}{2}-1} 
\frac{(1-v)^{\frac{\mu}{2}-b}Ê}{\lambda^{2b+\mu}(1+v\lambda^{-2})^{b+\frac{\mu}{2}+\frac{1}{2}}}\,.
$$
To derive the desired formula, we finally use a classical integral representation of 
${}_2F_1$, together with
$${}_2F_1\left(\alpha, \beta,\gamma;-z\right) = (1+z)^{\gamma-\alpha-\beta} {}_2F_1(\gamma-\alpha, \gamma-\beta, \gamma; -z)\,;$$
see formula (9.5.3) in Lebedev \cite{Leb}. We leave the details to the reader.


\end{document}